\documentclass[graybox]{svmult}

\usepackage{mathptmx}       
\usepackage{helvet}         
\usepackage{courier}        
\usepackage{makeidx}         
\usepackage{graphicx}        
\usepackage{multicol}        
\usepackage{epsfig}
\makeindex             

\begin{document}

\title{Showcase of Blue Sky Catastrophes}

\author{Leonid Shilnikov, Andrey Shilnikov and Dmitry Turaev}
\institute{Leonid Shilnikov \at Institute for Applied Mathematics \& Cybernetics,
10 Ulyanov Street, Nizhny Novgorod, 603005, Russia\\
 Andrey Shilnikov \at Neuroscience Institute, and
Department of Mathematics and Statistics,
Georgia State University, Atlanta,  100 Piedmont Ave SE, Atlanta, GA, 30303, USA, \email{ashilnikov@gsu.edu}, \\
Dmitry Turaev \at Department of Mathematics,
Imperial College, London, SW7 2AZ, UK, \email{d.turaev@imperial.ac.uk}}
%
%
\maketitle

\abstract{Let a system of differential equations possess a saddle-node periodic orbit such that
every orbit in its unstable manifold is homoclinic, i.e. the unstable manifold is a subset of the (global)
stable manifold. We study several bifurcation cases where the splitting of such a homoclinic connection
causes the Blue Sky Catastrophe, including the onset of complex dynamics.
The birth of an invariant torus or a Klein bottle is also described.}

\section{Introduction}

In the pioneering works by A.A.~Andronov and E.A.~Leontovich \cite{AL1,AL2}
all main bifurcations of stable periodic orbits of dynamical systems in a plane had been studied: the emergence
of a limit cycle from a weak focus,
the saddle-node bifurcation through a merger of a stable limit cycle with an unstable one and their consecutive
annihilation, the birth of a limit cycle from a separatrix loop to a saddle, as well as from a separatrix loop
to a saddle-node equilibrium. Later, in the 50-60s these bifurcations were generalized for the multi-dimensional case,
along with two additional bifurcations: period doubling and the birth of a two-dimensional torus. Apart from that,
in \cite{lp1,lp2} L.~Shilnikov had studied the main bifurcations of saddle periodic orbits out of homoclinic loops
to a saddle and discovered a novel bifurcation of homoclinic loops to a saddle-saddle\footnote{an equilibrium state, alternatively called a Shilnikov saddle-node,
due to a merger of two saddles of different topological types}.

Nevertheless, an open problem still remained: could there be other types of codimension-one bifurcations
of periodic orbits? Clearly, the emphasis was put on bifurcations of {\em stable} periodic orbits, as only
they generate robust self-sustained periodic oscillations, the original paradigm of
nonlinear dynamics. One can pose the problem as follows:\\ {\em
In a one-parameter family $X_{\mu}$ of systems of differential equations,
can both the period and the length of a structurally stable periodic orbit ${\cal L}_\mu$
tend to infinity as the parameter $\mu$ approaches some bifurcation value, say $\mu_0=0$?} \\
Here, structural stability means that none of the multipliers of the periodic orbit ${\cal L}_\mu$ crosses the unit circle, i.e. ${\cal L}_\mu$
does not bifurcate at $\mu\neq\mu_0$. Of particular interest is the case where ${\cal L}_\mu$ is stable, i.e.
all the multipliers are strictly inside the unit circle.

A similar formulation was given by J.~Palis and Ch.~Pugh \cite{PP} (notable Problem~37), however the structural stability
requirement was missing there. Exemplary bifurcations of a periodic orbit whose period becomes arbitrarily large
while the length remains finite as the bifurcation moment is approached are
a homoclinic bifurcation of a saddle with a negative saddle value and that
of a saddle-node \cite{lp0,book2}. These were well-known at the time, so in \cite{PP} an additional condition
was imposed, in order to ensure that the sought bifurcation is really of a new type: the periodic orbit ${\cal L}_\mu$ must stay
away from any equilibrium states (this would immediately imply that the length of the orbit grows to infinity
in proportion to the period). As R.~Abraham put it, the periodic orbit must ``disappear in the blue-sky'' \cite{Ab}.

In fact, a positive answer to ``Problem 37'' could be found in an earlier paper \cite{F}. In explicit form, a solution
was proposed by  V.~Medvedev \cite{Me}. He constructed examples of flows on a torus and a Klein bottle with
stable limit cycles whose lengths and periods tend to infinity as $\mu\to\mu_0$, while at $\mu=\mu_0$
both the periodic orbits disappear and new, structurally unstable saddle-node periodic orbits appear
(at least two of them, if the flow is on a torus). The third example of \cite{Me} was a flow on a 3-dimensional torus
whose all orbits are periodic and degenerate, and for the limit system the torus is foliated by two-dimensional invariant tori.

Medvedev's examples are not of codimension-1: this is obvious for the torus case that requires at least two saddle-nodes, i.e.
$X_{\mu_0}$ is of codimension 2 at least. In case of the Klein bottle one may show \cite{book2,AfS,TSh3,Li,Il}
that for a generic perturbation of the Medvedev family the periodic orbits existing at $\mu\neq\mu_0$
will not remain stable for all $\mu$ as they undergo an infinite sequence of
forward and backward period-doubling bifurcations (this is a typical  behavior of fixed points of a non-orientable
diffeomorphism of a circle).

\begin{figure}[htb!]
\begin{center}
\includegraphics[width=0.8\textwidth]{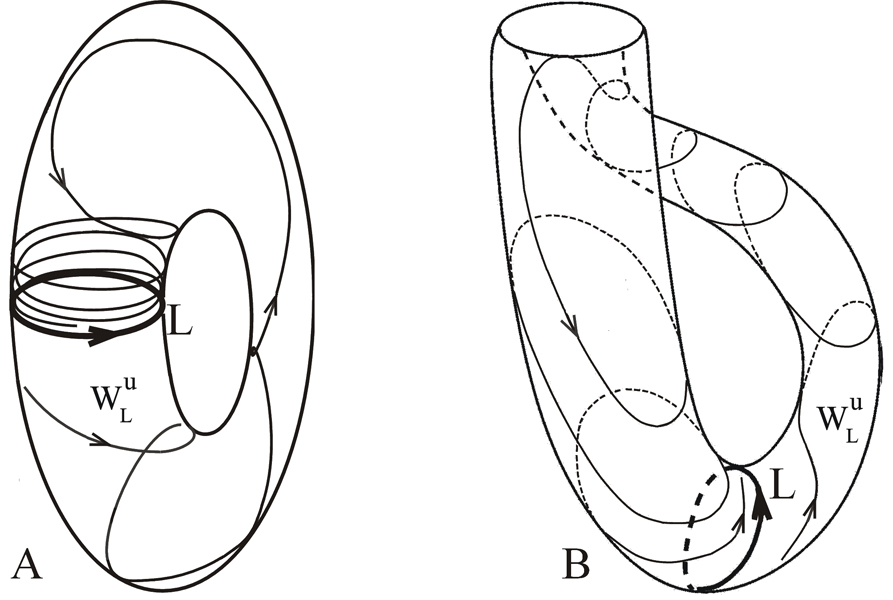}
\end{center}
\caption{Two cases of the unstable manifold $W^u_L$ homoclinic to the saddle-node periodic orbit $L$:
a 2D torus (A) or a Klein bottle (B).}
\label{fig1}
\end{figure}

A blue-sky catastrophe of codimension 1 was found only in 1995 by L.~Shilnikov and D.~Turaev \cite{TSh3,TSh1,TSh2,ShT}.
The solution was based on the study of bifurcations of a saddle-node periodic orbit whose entire unstable manifold
is homoclinic to it. The study of this bifurcation was initiated by V.~Afraimovich and L.~Shilnikov \cite{AfS,AfS1,AfS2,AfS3}
for the case where the unstable manifold of the saddle-node is  a torus or a Klein bottle (see Fig.~\ref{fig1}).
As soon as the saddle-node disappears, the Klein bottle may persist, or it may break down to cause chaotic
dynamics  in the system \cite{AfS4,NPT,TSh,Sync}. In these works, most of attention was paid to the torus case,
as its breakdown provides a geometrical model of the quasiperiodicity-toward-chaos transition encountered
universally in Nonlinear Dynamics, including the onset of turbulence \cite{Sh00}.

\begin{figure}[htb!]
\begin{center}
\includegraphics[height=0.45\textheight]{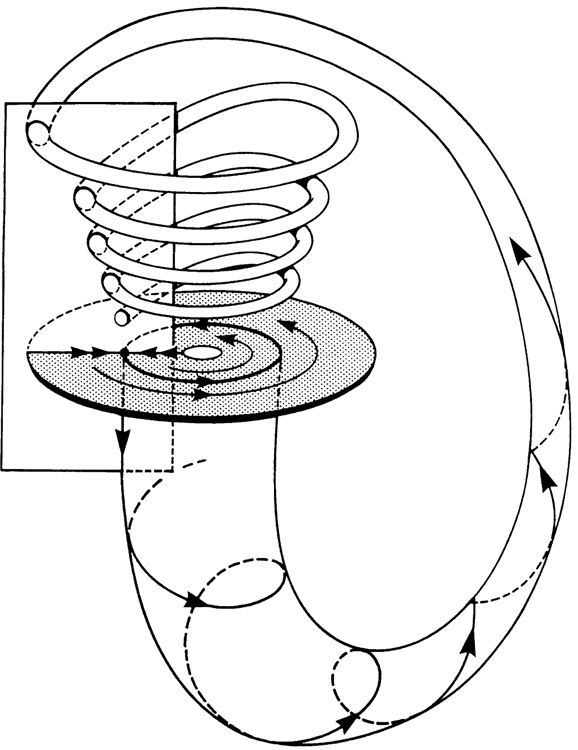}
\end{center}
\caption{Original construction of the blue sky catastrophe from \cite{TSh1}.}
\label{fig2}
\end{figure}

In the hunt for the blue sky catastrophe, other distinct configurations of the unstable manifold of the saddle-node were suggested
in \cite{TSh1}. In particular, it was shown that in the phase space of dimension 3 and higher
the homoclinic trajectories may spiral back onto the saddle-node orbit in the way shown in Fig.~\ref{fig2}.
If we have a one-parameter family $X_\mu$ of systems of differential equations
with a saddle-node periodic orbit at $\mu=\mu_0$ which possesses this special kind of the homoclinic unstable
manifold and satisfy certain additional conditions, then as the saddle-node disappears the inheriting attractor
consists of a single stable periodic orbit ${\cal L}_\mu$ which
undergoes no bifurcation as $\mu\to\mu_0$ while its length tends to infinity. Its topological limit, $M_0$, is
the entire unstable manifold of the saddle-node periodic orbit.

The conditions found in \cite{TSh1} for the behavior of the homoclinic orbits ensuring the blue-sky catastrophe are open,
i.e. a small perturbation of the one-parameter family $X_\mu$ does not destroy the construction. This implies
that such a  blue-sky catastrophe occurs any time a family of systems of differential equations crosses
the corresponding codimension-1 surface in the Banach space of smooth dynamical systems. This surface constitutes
a stability boundary for periodic orbits. This boundary is drastically new comparable to those
known since the 30-60s and has no analogues in planar systems. There are reasons to conjecture that this
type of the blue-sky catastrophe closes the list of main stability boundaries for periodic orbits (i.e. any new stability
boundary will be of codimension higher than 1).

In addition, another version of blue-sky catastrophe leading to the birth of a uniformly-hyperbolic strange attractor
(the Smale-Williams solenoid \cite{Sm,W}) was also discovered in \cite{TSh1,TSh2}. This codimension-1 bifurcation
of a saddle-node which corresponds yet to a different configuration of the homoclinic unstable manifold of the periodic orbit
(the full classification is presented in \cite{book2}). Here, the structurally stable attractor existing all the way
up to $\mu=\mu_0$ does not bifurcate so that the length of each and every (saddle) periodic orbit in it tends to
infinity as $\mu\to\mu_0$.

Initially we believed that the corresponding configuration of the unstable manifold would be too exotic for
the blue-sky catastrophe to occur naturally in a plausible system. In contrast, soon after, a first explicit
example of the codimension-1 blue-sky catastrophe
was proposed by N.~Gavrilov and A.~Shilnikov \cite{GSh}, in the form of a family of 3D systems of differential equations
with polynomial right-hand sides. A real breakthrough came in when the blue-sky catastrophe has turned out to be a typical
phenomenon for slow-fast systems. Namely, in \cite{book2,mmj} we described a number of very general scenarios leading
to the blue-sky catastrophe in such
systems with at least two fast variables; for systems with one fast variable the blue-sky catastrophe was found in \cite{GKR}.
In this way, the blue-sky catastrophe has found numerous applications in mathematical neuroscience, namely, it explains a smooth
and reversible transition between tonic spiking and bursting in exact Hodgkin-Huxley type models of
interneurons \cite{leech1,leech2} and in mathematical models of square-wave bursters \cite{hr}.
The great variability of the burst duration near the blue-sky catastrophe
was shown to be the key mechanism ensuring the diversity of rhythmic patterns generated by small neuron complexes
that control invertebrate locomotion \cite{DG1,DG2,DG3}.

In fact, the term ``blue sky catastrophe" should be naturally treated in a broader way. Namely, under this term we allow
to embrace a whole class of dynamical phenomena that all are due to the existence of a stable (or, more generally, structurally stable)
periodic orbit, ${\cal L}_\mu$, depending continuously on the parameter $\mu$ so that both, the length and the period of ${\cal L}_\mu$ tend
to infinity as the bifurcation parameter value is reached. As for the topological limit, $M_0$,
of the orbit ${\cal L}_\mu$ is concerned, it may possess a rather degenerate structure that does not prohibit $M_0$ from
having equilibrium states included. As such, the periodic regime ${\cal L}_\mu$ could emerge as a composite construction
made transiently of several quasi-stationary states: nearly constant, periodic, quasiperiodic, and even chaotic fragments.
As one of the motivations (which we do not pursue here) one may think on slow-fast model where the fast
3D dynamics is driven by a periodic motion in a slow subsystem.

\section{Results}

In this paper we focus on an infinitely degenerate case where $M_0$
is comprised of a saddle periodic orbit with a continuum of homoclinic trajectories.
Namely, we consider a one-parameter family of sufficiently smooth systems of differential equations
$X_\mu$ defined in $R^{n+1}$, $n\geq 2$, for which we need to make a number of assumptions as follows.\\

\noindent {\bf (A)} There exists a saddle periodic orbit $L$ (we assume the period equals $2\pi$) with the
multipliers\footnote{the eigenvalues of the linearization of the Poincare map} $\rho_1,\dots,\rho_n$. Let the multipliers satisfy
\begin{equation}\label{rh1}
 \max_{i=2,\dots,n-1} |\rho_i|<|\rho_1|<\;1\;<|\rho_n|.
\end{equation}
Once this property is fulfilled at $\mu=0$, it implies that the saddle periodic orbit $L=L_\mu$ exists
for all small $\mu$ and smoothly depends on $\mu$. Condition (\ref{rh1}) also holds for all small $\mu$.
This condition implies that the
stable manifold $W^s_\mu$ is $n$-dimensional\footnote{the intersection of $W^s_\mu$ with any cross-section to $L_\mu$ is $(n-1)$-dimensional}
and the unstable manifold $W^u_\mu$ is two-dimensional.
If the unstable multiplier $\rho_n$ is positive (i.e. $\rho_n>1$), then
the orbit $L_\mu$ divides $W^u_\mu$ into two halves, $W^+_\mu$ and $W^-_\mu$, so
$W^u_\mu=L_\mu\cup W^+_\mu\cup W^-_\mu$. If $\rho_n$ is negative ($\rho_n<-1$), then
$W^u_\mu$ is a M\"obius strip, so $L_\mu$ does not divide $W^u_\mu$; in this case we denote
$W^+_\mu=W^u_\mu\backslash L_\mu$.

Concerning the stable manifold, condition (\ref{rh1}) implies that in $W^s_\mu$ there
exists (at $n\geq 3$) an $(n-1)$-dimensional strong-stable invariant manifold $W^{ss}_\mu$ whose
tangent at the points of $L_\mu$ contains the eigen-directions
corresponding to the multipliers $\rho_2,\dots,\rho_{n-1}$, and the orbits in $W^s_\mu\backslash W^{ss}_\mu$
tend to $L_\mu$ along the direction which correspond to the leading multiplier $\rho_1$.\\

\noindent {\bf (B)} At $\mu=0$ we have $W^+_0\subset W^s_0\backslash W^{ss}_0$,
i.e. we assume that {\em all} orbits from $W^+_0$ are
homoclinic to $L$. Moreover, as $t\to +\infty$, they tend to $L$ along the leading direction.\\

\noindent {\bf (C)} We assume that the flow near $L$ contracts three-dimensional volumes, i.e.
\begin{equation}\label{contr}
|\rho_1\rho_n| <1.
\end{equation}
This condition is crucial, as the objects that we obtain by bifurcations of the homoclinic surface $W^+_0\cup L$ are meant to
be attractors. Note that this condition is similar to the negativity of the saddle value condition from the theory of
homoclinic loops to a saddle equilibrium \cite{AL1,AL2,lp0}, see (\ref{sadl}).\\

\noindent {\bf (D)} We assume that one can introduce linearizing coordinates near $L$. Namely, a small neighborhood $U$ of $L$
is a solid torus homeomorphic to $S^1\times R^n$, i.e. we can coordinatize it by an angular variable $\theta$
and by normal coordinates $u\in R^n$. Our assumption is that these coordinates are chosen so that
the system in the small neighborhood of $L$ takes the form
\begin{equation}\label{lfr}
 \dot u=C(\theta,\mu) u, \qquad \dot \theta=1,
\end{equation}
where $C$ is $2\pi$-periodic in $\theta$. The smooth linearization is not always possible, and our results
can be obtained without this assumption. We, however, will avoid discussing the general case here, in order
to make the construction more transparent.

It is well-known that by a $4\pi$-periodic transformation of the coordinates $u$ system (\ref{lfr}) can be brought to
the time-independent form. Namely, we may write the system as follows
\begin{equation}\label{lcfr}
\begin{array}{l}
\dot x=-\lambda(\mu) x, \qquad \dot y=B(\mu) y,\\
\dot z=\gamma(\mu) z,\\
\dot \theta=1,\end{array}
\end{equation}
where $x\in R^1$, $y\in R^{n-2}$, $z\in R^1$, and $\lambda=-\frac{1}{2\pi}\ln|\rho_1|>0$, $\gamma=\frac{1}{2\pi}\ln|\rho_n|>0$
and, if $n\geq 2$, $B(\mu)$ is an $(n-2)\times(n-2)$-matrix such that
\begin{equation}\label{nev}
\|e^{Bt}\|=o(e^{-\lambda t}) \qquad (t\to+\infty).
\end{equation}
Note also that condition {\bf (C)} implies
\begin{equation}\label{sadl}
\gamma-\lambda<0.
\end{equation}
By (\ref{lcfr}), the periodic orbit $L(\mu)$ is given by $x=0$, $y=0$, $z=0$, its local stable manifold is given
by $z=0$, and the leading direction in the stable manifold is given by $y=0$; the local unstable manifold
is given by $\{x=0,y=0\}$.

Recall that the $4\pi$-periodic transformation we used to bring system (\ref{lfr}) to the autonomous form (\ref{lcfr})
is, in fact, $2\pi$-periodic or $2\pi$-antiperiodic. Namely, the points $(\theta,x,z,y)$ and $(\theta+2\pi,\sigma(x,z,y))$
are equal (they represent the same point in the solid torus $U$), where $\sigma$
is an involution which changes signs
of some of the coordinates $x,z,y_1,\dots,y_{n-2}$. More precisely, $\sigma$ changes the orientation of each of the directions
which correspond to the real negative multipliers $\rho$. In particular, if all the multipliers $\rho$ are positive, then $\sigma$
is the identity, i.e. our coordinates are $2\pi$-periodic in this case.\\

\begin{figure}[htb!]
\begin{center}
\includegraphics[width=0.8\textwidth]{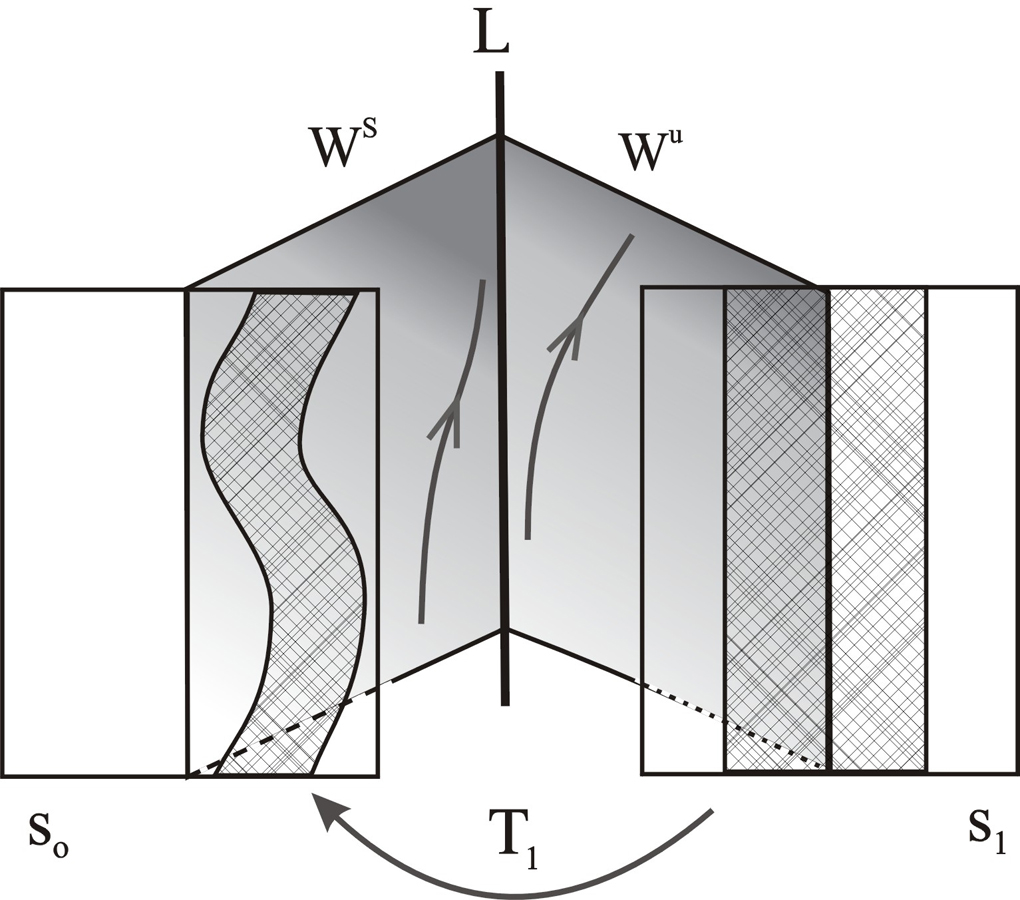}
\end{center}
\caption{Poincar\'e map $T_1$ takes a cross-section $S_1$ transverse to the unstable manifold $W^u$
to a cross-section $S_0$ transverse to the stable manifold $W^s$.}
\label{fig3}
\end{figure}

\noindent {\bf (E)} Consider two cross-sections $S_0:\{x=d,\quad \|y\|\leq \varepsilon_1,\quad |z|\leq \varepsilon_1\}$ and
$S_1:\{z=d,\quad \|y\|\leq\varepsilon_2,\quad |x|\leq\varepsilon_2\}$ for some small positive $d$ and $\varepsilon_{1,2}$.
Denote the coordinates on $S_0$ as $(y_0,z_0,\theta_0)$ and the coordinates on $S_1$ as $(x_1,y_1,\theta_1)$.
The set $S_0$ is divided by the stable manifold $W^s$ into two regions, $S_0^+:\{z_0>0\}$ and $S_0^-:\{z_0<0\}$.
Since $W^+_0\subset W^s_0$ by assumption 2, it follows that the orbits starting at $S_1$
define a smooth map $T_1:S_1\to S_0$ (see Fig.~\ref{fig3}) for all small $\mu$:
\begin{equation}\label{glom}
\begin{array}{l}
z_0 =f(x_1,y_1,\theta_1,\mu)\\
y_0 =g(x_1,y_1,\theta_1,\mu)\\
\theta_0 =m\theta_1 + h(\theta_1,\mu)+\tilde h(x_1,y_1,\theta_1,\mu),
\end{array}
\end{equation}
where $f,g,h,\tilde h$ are smooth functions $4\pi$-periodic in $\theta_1$, and the function $\tilde h$ vanishes at $(x_1=0,y_1=0)$.
Condition $W^+_0\subset W^s_0$ reads as
$$f(0,0,\theta_1,0)\equiv 0.$$
We assume that
\begin{equation}\label{qqfff}
f(0,0,\theta_1,\mu)=\mu\alpha(\theta_1,\mu),
\end{equation}
where
\begin{equation}\label{alpt}
\alpha(\theta_1,\mu)>0
\end{equation}
for all $\theta_1$, i.e. {\em all the homoclinics are split simultaneously and in the same direction}, and
the intersection $W^+_\mu\cap S_0$ moves inside $S_0^+$ with a non-zero velocity as $\mu$ grows across zero.

The coefficient $m$ in the last equation of (\ref{glom}) is an integer. In order to see this, recall that
two points $(\theta,x,z,y)$ and $(\hat\theta,\hat x,\hat z,\hat y)$ in $U$ are the same if and only if
$\hat\theta=\theta+2\pi k, (\hat x,\hat z,\hat y)=\sigma^k (x,z,y)$ for an integer $k$. Thus, if we increase $\theta_1$ to $4\pi$
in the right-hand side of (\ref{glom}), then the corresponding value of $\theta_0$ in the left-hand side
may change only to an integer multiple of $2\pi$, i.e. $m$ must be an integer or a half-integer. Let us show that
the half-integer $m$ are forbidden by our assumption (\ref{alpt}). Indeed, if the multiplier $\rho_n$ is positive, then
the involution $\sigma$ keeps the corresponding variable $z$ constant. Thus, $(z=d,\theta=\theta_1, x=0, y=0)$ and
$(z=d,\theta=\theta_1+2\pi, x=0, y=0)$ correspond, in this case, to the same point on $W^+_\mu\cap S_1$, hence their image
by (\ref{glom}) must give the same point on $S_0$, i.e. the corresponding values of $\theta_0$ must differ on an integer multiple
of $2\pi$, which means that $m$ must be an integer. If $\rho_n<0$, then $\sigma$ changes the sign of $z$, i.e. if two values of
$\theta_0$ which correspond to the same point on $S_0$ differ on $2\pi k$, the corresponding values of $z$ differ to a factor of
$(-1)^k$. Now, since the increase of $\theta_1$ to $4\pi$ leads to the increase of $\theta_0$ to $4\pi m$ in (\ref{glom}),
we find that $f(0,0,4\pi,\mu)=(-1)^{2m}f(0,0,0,\mu)$ in the case $\rho_n<0$. This implies
that if $m$ is a half-integer, then $f(0,0,\theta)$ must have zeros at any $\mu$ and (\ref{alpt}) cannot be satisfied.

The number $m$ determines the shape of $W^+\cap S_0$. Namely, the equation of the curve $W^+_0\cap S_0$ is
$$\theta_0 =m\theta_1 + h_1(\theta_1,0),\qquad y_0 =g(0,0,\theta_1,0), \qquad z_0=0,$$
so $|m|$ defines the homotopic type of this curve in $S_0\cap W^s_0$, and the sign of $m$ is responsible for the orientation.
In the case $n=2$, i.e. when the system is defined in $R^3$, the only possible case is $m=1$. At $n=3$ (the system in $R^4$)
the curve $W^+_0\cap S_0$ lies in the two-dimensional intersection of $W^s$ with $S_0$. This is either an annulus (if $\rho_1>0$),
or a M\"obius strip (if $\rho_1<0$). Since the smooth curve $W^+_0\cap S_0$ cannot have self-intersections, it follows that
the only possible cases are $m=0,\pm1$ when $W^s\cap S_0$ is a two-dimensional annulus and $m=0,\pm1,\pm2$ when
$W^+_0\cap S_0$ is a M\"obius strip. At large $n$ (the system in $R^5$ and higher) all integer values of $m$ are possible.\\~\\

Now we can formulate the main results of the paper.\\

\noindent{\bf Theorem.}
Let conditions {\bf (A-E)} hold. Consider a sufficiently small neighborhood $V$ of the homoclinic surface $\Gamma=W^+_0\cap L$.\\
\begin{enumerate}
\item If $m=0$ and, for all $\theta$,
\begin{equation}\label{bsky}
|h'(\theta,0)-\frac{\alpha'(\theta,0)}{\gamma \alpha(\theta,0)}|<1,
\end{equation}
then a single stable periodic orbit ${\cal L}_\mu$ is born as $\Gamma$ splits. The orbit ${\cal L}_\mu$ exists at all small $\mu>0$; its
period and length tend to infinity as $\mu\to+0$.
All orbits which stay in $V$ for all positive times and
which do not lie in the stable manifold of the saddle orbit $L_\mu$ tend to ${\cal L}_\mu$.\\
\item If $|m|=1$ and, for all $\theta$,
\begin{equation}\label{tor}
1+m \left[h'(\theta,0)-\frac{\alpha'(\theta,0)}{\gamma \alpha(\theta,0)}\right]>0,
\end{equation}
then a stable two-dimensional invariant torus (at $m=1$) or a Klein bottle (at $m=-1$) is born as $\Gamma$ splits. It exists
at all small $\mu>0$ and attracts all the orbits which stay in $V$ and which do not lie in the stable manifold of $L_\mu$.\\
\item If $|m|\geq 2$ and, for all $\theta$,
\begin{equation}\label{hypat}
|m+h'(\theta,0)-\frac{\alpha'(\theta,0)}{\gamma \alpha(\theta,0)}|>1,
\end{equation}
then, for all small $\mu>0$, the system has a hyperbolic attractor (a Smale-Williams solenoid) which is an $\omega$-limit set
for all orbits which stay in $V$ and
which do not lie in the stable manifold of $L_\mu$. The flow on the attractor is topologically conjugate to
suspension over the inverse spectrum limit of a degree-$m$ expanding map of a circle. At $\mu=0$, the attractor
degenerates into the homoclinic surface $\Gamma$.\\
\end{enumerate}

\noindent{\em Proof.} Solution of (\ref{lcfr}) with the initial conditions $(x_0=d,y_0,z_0,\theta_0)\in S_0$ gives
$$\begin{array}{l}
x(t)=e^{-\lambda t} d, \qquad y(t)=e^{B t} y_0,\\
z(t)=e^{\gamma t} z_0,\\
\theta(t)=\theta_0+t.\end{array}
$$
The flight time to $S_1$ is found from the condition
$$d = e^{\gamma t} z_0,$$
which gives $\displaystyle t=-\frac{1}{\gamma}\ln\frac{z_0}{d}$. Thus the orbits in $U$ define the map $T_0: S_0^+\to S_1$:
$$\begin{array}{l}
x_1=d^{1-\nu} z_0^\nu, \qquad y_1=Q(z_0) y_0,\\
\theta_1=\theta_0-\frac{1}{\gamma}\ln\frac{z_0}{d}\end{array}
$$
where $\nu=\lambda/\gamma>1$ and $\|Q(z_0)\|=o(z_0^\nu)$ (see (\ref{nev}),(\ref{sadl})). By (\ref{glom}), we may write the
map $T=T_0T_1$ on $S_1$ as follows (we drop the index ``$1$''):
$$\begin{array}{l}
\bar x=d^{1-\nu} (\mu\alpha(\theta,\mu)+O(x,y))^\nu, \qquad \bar y=Q(\mu\alpha+O(x,y)) g(x,y,\theta,\mu),\\
\bar\theta=m\theta+h(\theta,\mu)-\frac{1}{\gamma}\ln(\frac{\mu}{d}\alpha(\theta,\mu)+O(x,y))+O(x,y).\end{array}
$$
For every orbit which stays in $V$, its consecutive intersections with the cross-section $S_1$ constitute an orbit of
the diffeomorphism $T$. Since $\nu>1$, the map $T$ is contracting in $x$ and $y$, and it is easy to see
that all the orbits eventually enter a neighborhood of $(x,y)=0$ of size $O(\mu^\nu)$.
We therefore rescale the coordinates $x$ and $y$ as follows:
$$x=d^{1-\nu}\mu^\nu X,\qquad y=\mu^\nu Y.$$
The map $T$ takes the form
\begin{equation}\label{mapt}
\begin{array}{l}
\bar X= \alpha(\theta,0)^\nu +o(1), \qquad \bar Y=o(1),\\
\bar\theta=\omega(\mu)+m\theta +h(\theta,0)-\frac{1}{\gamma}\ln\alpha(\theta,0)+o(1),
\end{array}
\end{equation}
where $o(1)$ stands for terms which tend to zero as $\mu\to+0$, along with their first derivatives,
and $\omega(\mu)=\frac{1}{\gamma}\ln(\mu/d)\to\infty$ as $\mu\to+0$. Recall that $\alpha>0$ for all $\theta$
and that $\alpha$ and $h$ are periodic in $\theta$.

\begin{figure}[htb!]
\begin{center}
\includegraphics[width=0.8\textwidth]{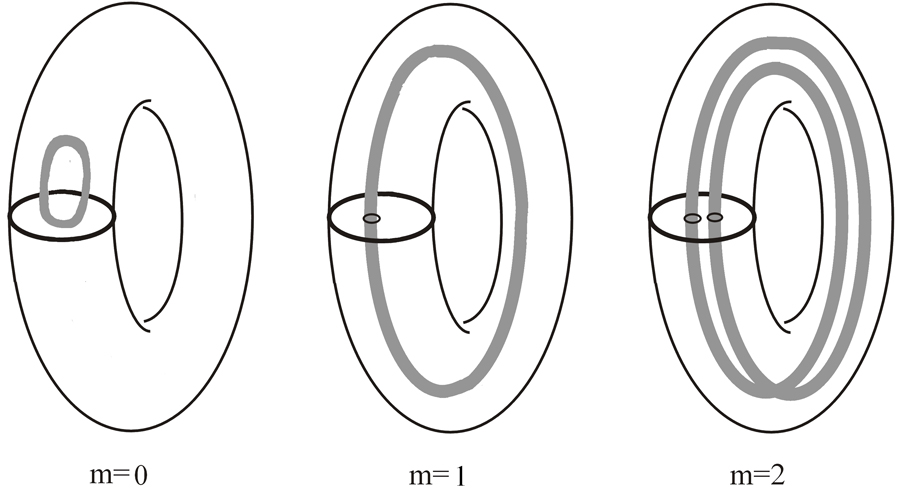}
\end{center}
\caption{Case $m=0$: the image of the solid torus is contractible to a point; case $m = 1$: contraction transverse to the longitude;
case $m = 2$: the solid-torus is squeezed,
doubly stretched and twisted within the original and so on, producing the solenoid in the limit.}
\label{fig4}
\end{figure}

It is immediately seen from (\ref{mapt}) that all orbits
eventually enter an invariant solid torus $\{|x-\alpha(\theta,0)^\nu|< K_\mu,\;\|y\|<K_\mu\}$
for appropriately chosen $K_\mu$, $K_\mu\to 0$ as $\mu\to +0$ (see Fig.~\ref{fig4}). Thus, there is an attractor in $V$ for
all small
positive $\mu$, and it merges into $\Gamma$ as $\mu\to+0$. Our theorem claims that the structure of the attractor depends
on the value of $m$, so we now consider different cases separately.

If $m=0$ and (\ref{bsky}) holds, then map (\ref{mapt}) is, obviously, contracting at small $\mu$,
hence it has a single stable fixed point. This fixed point corresponds to the sought periodic orbit
$A_\mu$. Its period tends to infinity as $\mu\to+0$: the orbit intersects both the
cross-sections $S_0$ and $S_1$, and the flight time from $S_0$ to $S_1$ is of order $\frac{1}{\gamma}|\ln\mu|$.
The length of the orbit also tends to infinity, since the phase velocity never vanishes in $V$.

In the case $m=\pm 1$ we prove the theorem by referring to the ``annulus principle'' of \cite{AfS3}. Namely, consider a map
$$\bar r=p(r,\theta),\qquad \bar\theta=q(r,\theta)$$
of a solid torus into itself (here $\theta$ is the angular variable and $r$ is the vector of normal variables).
Let the map $r\mapsto p(r,\theta)$ be a contraction for every fixed $\theta$, i.e.
$$\left\|\frac{\partial p}{\partial r}\right\|_\circ<1$$
(where by $\|\cdot\|_\circ$ we denote the supremum of the norm over the solid torus under consideration)
and let the map $\theta\mapsto q(r,\theta)$ be a diffeomorphism of a circle for every fixed $r$. Then
it is well-known \cite{AfS3,book2} that if
$$1-\left\|\left(\frac{\partial q}{\partial \theta}\right)^{-1}\right\|_\circ \cdot
\left\|\frac{\partial p}{\partial r}\right\|_\circ >
2\sqrt{\left\|\left(\frac{\partial q}{\partial \theta}\right)^{-1}\right\|_\circ \cdot
\left\|\frac{\partial q}{\partial r}\right\|_\circ
\left\|\frac{\partial p}{\partial \theta}\left(\frac{\partial q}{\partial \theta}\right)^{-1}\right\|_\circ},$$
then the map has a stable, smooth, closed invariant curve $r=r^*(\theta)$ which attracts all orbits from the solid torus.
These conditions are clearly satisfied by map (\ref{mapt}) at $|m|=1$ if (\ref{tor}) is true (here $r=(X,Y)$,
$p=(\alpha(\theta,0)^\nu +o(1), o(1))$, $q=\omega(\mu)+m\theta +h(\theta,0)-\frac{1}{\gamma}\ln\alpha(\theta,0)+o(1)$).
Thus, the map $T$ has a a closed invariant curve in this case. The restriction of $T$ to the invariant curve preserves
orientation if $m=1$, while at $m=-1$ it is orientation-reversing. Therefore, this invariant curve on the cross-section
corresponds to an invariant torus of the flow at $m=1$ or to a Klein bottle at $m=-1$.

It remains to prove the theorem for the case $|m|\geq 2$. The proof is based on the following result.\\
\noindent{\bf Lemma.} Consider a diffeomorphism $T:(r,\theta)\mapsto (\bar r,\bar\theta)$ of a solid torus, where
\begin{equation}\label{maptr}
\bar r=p(r,\theta),\qquad \bar\theta=m\theta+s(r,\theta)=q(r,\theta),
\end{equation}
where $s$ and $p$ are periodic functions of $\theta$
Let $|m|\geq 2$, and
\begin{equation}\label{frc}
\left\|\frac{\partial p}{\partial r}\right\|_\circ <1,
\end{equation}
\begin{equation}\label{cndir}
\left(1-\left\|\frac{\partial p}{\partial r}\right\|_\circ\right)
\left(1-\left\|\left(\frac{\partial q}{\partial \theta}\right)^{-1}\right\|_\circ\right)>
\left\|\frac{\partial p}{\partial \theta}\right\|_\circ\; \left\|\left(\frac{\partial q}{\partial \theta}\right)^{-1}
\frac{\partial q}{\partial r}\right\|_\circ.
\end{equation}
Then the map has a uniformly-hyperbolic attractor, a Smale-Williams solenoid,
on which it is topologically conjugate to the inverse spectrum limit
of $\bar \theta=m\theta$, a degree-$m$ expanding map of the circle.\\

\noindent{\em Proof.} It follows from (\ref{frc}),(\ref{cndir}) that
$\|(\frac{\partial q}{\partial \theta})^{-1}\|$ is uniformly bounded. Therefore, $\theta$
is a uniquely defined smooth function of $(\bar\theta, r)$, so we may rewrite (\ref{maptr})
in the ``cross-form''
\begin{equation}\label{crmps}
\bar r=p^\times(r,\bar\theta),\qquad \theta=q^\times(r,\bar\theta),
\end{equation}
where $p^\times$ and $q^\times$ are smooth functions. It is easy to see that conditions (\ref{frc}),
(\ref{cndir}) imply
\begin{equation}\label{frc0}
\left\|\frac{\partial p^\times}{\partial r}\right\|_\circ <1,\qquad
\left\|\frac{\partial q^\times}{\partial \bar\theta}\right\|_\circ <1
\end{equation}
\begin{equation}\label{cncrs}
\left(1-\left\|\frac{\partial p^\times}{\partial r}\right\|_\circ\right)
\left(1-\left\|\frac{\partial q^\times}{\partial \theta}\right\|_\circ\right)\geq
\left\|\frac{\partial p^\times}{\partial \bar\theta}\right\|_\circ\; \left\|\frac{\partial q^\times}{\partial r}\right\|_\circ.
\end{equation}
These inequalities imply the uniform hyperbolicity of the map $T$ (note that (\ref{cndir}) coincides with
the hyperbolicity condition for the Poincare map for the Lorenz attractor from \cite{ABS}).
Indeed, it is enough to show that there exists
$L>0$ such that the derivative $T'$ of $T$ takes every cone $\|\Delta r\|\leq L\|\Delta \theta\|$ inside
$\|\Delta \bar r\|\leq L\|\Delta \bar \theta\|$ and is uniformly expanding in $\theta$ in this cone,
and that the inverse of $T'$ takes
every cone $\|\Delta \bar\theta\|\leq L^{-1}\|\Delta \bar r\|$ inside
$\|\Delta \theta\|\leq L^{-1}\|\Delta r\|$ and is uniformly expanding in $r$ in this cone.
Let us check these properties. When $\|\Delta r\|\leq L\|\Delta \theta\|$, we find from (\ref{crmps}) that
\begin{equation}
 \|\Delta\theta\|\leq \frac{\left\|\frac{\partial q^\times}{\partial \bar\theta}\right\|_\circ}
{1-L\left\|\frac{\partial q^\times}{\partial r}\right\|_\circ} \|\Delta\bar\theta\|
\end{equation}
and
\begin{equation}
 \|\Delta\bar r\|\leq \left\{\frac{L \left\|\frac{\partial p^\times}{\partial r}\right\|_\circ
 \left\|\frac{\partial q^\times}{\partial \bar\theta}\right\|_\circ}
{1-L\left\|\frac{\partial q^\times}{\partial r}\right\|_\circ} +
\left\|\frac{\partial p^\times}{\partial \bar\theta}\right\|_\circ\right\}
\|\Delta\bar\theta\|.
\end{equation}
Similarly, if
$\|\Delta \bar\theta\|\leq L^{-1}\|\Delta \bar r\|$, we find from (\ref{crmps}) that
\begin{equation}
 \|\Delta\bar r\|\leq \frac{\left\|\frac{\partial p^\times}{\partial r}\right\|_\circ}
{1-L^{-1}\left\|\frac{\partial p^\times}{\partial \bar\theta}\right\|_\circ} \|\Delta r\|
\end{equation}
and
\begin{equation}
 \|\Delta\theta\|\leq \left\{\frac{L^{-1} \left\|\frac{\partial q^\times}{\partial \bar\theta}\right\|_\circ
\left\|\frac{\partial p^\times}{\partial r}\right\|_\circ}
{1-L^{-1}\left\|\frac{\partial p^\times}{\partial \bar\theta}\right\|_\circ} +
\left\|\frac{\partial q^\times}{\partial r}\right\|_\circ\right\}
\|\Delta r\|.
\end{equation}
Thus, we will prove hyperbolicity if we show that there exists $L$ such that
$$\left\|\frac{\partial q^\times}{\partial \bar\theta}\right\|_\circ < 1-
L\left\|\frac{\partial q^\times}{\partial r}\right\|_\circ$$
and
$$\left\|\frac{\partial p^\times}{\partial r}\right\|_\circ < 1-
L^{-1}\left\|\frac{\partial p^\times}{\partial \bar\theta}\right\|_\circ.$$
These conditions are solved by any $L$ such that
$$\frac{\left\|\frac{\partial p^\times}{\partial \bar\theta}\right\|_\circ}
{1-\left\|\frac{\partial p^\times}{\partial r}\right\|_\circ}<L<
\frac{1-\left\|\frac{\partial q^\times}{\partial \bar\theta}\right\|_\circ}
{\left\|\frac{\partial q^\times}{\partial r}\right\|_\circ}.$$
It remains to note that such $L$ exist indeed when (\ref{frc0}),(\ref{cncrs}) are satisfied.

We have proved that the attractor $A$ of the map $T$ is uniformly hyperbolic. Such attractors are structurally stable,
so $T|_A$ is topologically conjugate to the restriction to the attractor of any diffeomorphism which
can be obtained by a continuous deformation of the map $T$ without violation of conditions
(\ref{frc}) and (\ref{cndir}). An obvious example of such a diffeomorphism is given by the map
\begin{equation}\label{epd}
\bar r=p(\delta r,\theta),\qquad \bar\theta=q(\delta r,\theta)
\end{equation}
for any $0<\delta\leq 1$. Fix small $\delta>0$ and consider a family of maps
$$\bar r=p(\delta r,\theta),\qquad \bar\theta=q(\varepsilon r,\theta),$$
where $\varepsilon$ runs from $\delta$ to zero. When $\delta$ is sufficiently small, every map in this family is a diffeomorphism
(otherwise we would get that the curve $\{\bar r=p(0,\theta), \bar\theta= q(0,\theta)\}$ would have points of self-intersection,
which is impossible since this curve is the image of the circle $r=0$ by the diffeomorphism $T$), and each satisfies
inequalities (\ref{frc}),(\ref{cndir}). This family is a continuous deformation of map (\ref{epd}) to the map
\begin{equation}\label{skd}
\bar r=p(\delta r,\theta),\qquad \bar\theta=q(0,\theta)=m\theta+s(0,\theta).
\end{equation}
Thus, we find that $T|_A$ is topologically conjugate to the restriction of diffeomorphism (\ref{skd}) to its attractor.
It remains to note that map (\ref{skd}) is a skew-product map of the solid torus, which contracts along the fibers $\theta=const$
and, in the base, it is an expanding degree-$m$ map of a circle. By definition, the attractor of such map is the sought
Smale-Williams solenoid \cite{Sm,W}. This completes the proof of the lemma.

Now, in order to finish the proof of the theorem, just note that map (\ref{mapt}) satisfies the conditions of the Lemma
when (\ref{hypat}) is fulfilled.

\section*{Acknowledgment}

This work was supported by RFFI Grant No.~08-01-00083 and the Grant 11.G34.31.0039 of the Government of the Russian Federation
for state support of ``Scientific research conducted under supervision of leading scientists in
Russian educational institutions of higher professional education"  (to L.S);  NSF grant DMS-1009591,
MESRF ``Attracting leading scientists to Russian universities" project 14.740.11.0919 (to A.S) and
the Royal Society Grant "Homoclinic bifurcations" (to L.S. and D.T.)

\end{document}